\documentclass[a4paper,11pt]{article}
\usepackage[hidelinks]{hyperref}
\hypersetup{
	colorlinks=true,   
	linktoc=all,        
	linkcolor=red,     
	citecolor=blue,}
\usepackage{cite}
\usepackage{float}
\usepackage{tikz}
\usepackage{graphicx}
\usepackage{amssymb}
\usepackage{latexsym, bm}
\usepackage{multicol}
\usepackage{indentfirst}
\usepackage{amssymb,amsfonts}
\usepackage{amsmath}
\usepackage{setspace}
\newtheorem{theorem}{Theorem}[section]
\newtheorem{lemma}{Lemma}[section]
\newtheorem{claim}{Claim}[section]

\newtheorem{corollary}{Corollary}[section]
\newtheorem{proposition}{Proposition}[section]

\newtheorem{problem}{Problem}[section]

\newcommand{\ignore}[1]{}

\begin{document}

\begin{spacing}{1.05}
\title{Asymptotically optimal Ramsey goodness of sparse graphs versus\\ odd cycles and paths}
\date{}

\author{
Chunchao Fan\footnote{Center for Discrete Mathematics, Fuzhou University, Fuzhou, 350108 P.~R.~China. Email: {\tt 1807951575@qq.com}.} \;\; and\;\;
Qizhong Lin\footnote{Center for Discrete Mathematics, Fuzhou University,
Fuzhou, 350108 P.~R.~China. Email: {\tt linqizhong@fzu.edu.cn}. Supported in part  by National Key R\&D Program of China (Grant No. 2023YFA1010202) and NSFC (No.\ 12571361).}
}

\maketitle
\begin{abstract}

 A fundamental problem in graph Ramsey theory is to determine, for sparse graphs $G$ on $n$ vertices, the minimal $n$ such that $G$ is Ramsey-good for odd cycles $C_k$ and paths $P_k$. Burr, Erd\H{o}s, Faudree, Rousseau, and Schelp (Trans. AMS 1982) addressed this problem, establishing bounds requiring $n = \Omega(k^{10})$ for odd cycles and $n = \Omega(k^{12})$ for paths.
We settle the asymptotic version of this problem, proving that these bounds are essentially tight: $n = \Omega(k)$ suffices for odd cycles and $n = \Omega(k^2)$ (or $n = \Omega(k)$ under  additional conditions) for paths. Specifically, we prove:

\begin{itemize}
    \item[(1)] For odd cycles $C_k$ ($k\ge3$), we prove $r(G, C_k) = 2n-1$ for any connected $n$-vertex graph $G$ satisfying the relaxed conditions $n = \Omega(k)$ and $e(G) \le (1 + O(1/k^2)) n$. 
        
    \item[(2)] For paths $P_k$ ($k\ge2$), we prove $r(G, P_k) = \max\{ n + \lfloor k/2\rfloor - 1, n + k - 2 - \alpha' - \gamma \}$ for any connected $n$-vertex graph $G$  satisfying one of the following: 

\quad (i) $n = \Omega(k^2)$ and $e(G) \le (1 + O(1/k^2)) n$;
 
\quad (ii) $n = \Omega(k)$, $\delta(G)\ge2$, $\alpha'\geq k/2$, and $e(G) \le (1 + O(1/k)) n$.

In the above, $\alpha'$ is the independence number of an appropriate subgraph of $G$ and $\gamma=0$ if $k-1$ divides $n+k-3-\alpha'$, and $\gamma=1$ otherwise.

\end{itemize}

Consequently, our results unify and generalize classical theorems on odd cycles due to Bondy and Erd\H{o}s (1973), Faudree and Schelp (1974), and Rosta (1973), and on paths due to Gerencs\'er and Gy\'arf\'as (1967), Faudree, Lawrence, Parsons and Schelp (1974), and Parsons (1974). The proofs feature two key innovations: a novel reconstruction of the end-edge matching and an enhancement of Burr et al.'s dichotomy lemma.

\medskip

\end{abstract}

\section{Introduction}

For graphs $H_1$ and $H_2$, the Ramsey number $r(H_1,H_2)$ is the smallest positive integer $N$ such that any graph $G$ on $N$ vertices contains $H_1$ as a subgraph, or its complement contains $H_2$ as a subgraph. The existence of the Ramsey number $r(H_1,H_2)$ follows from \cite{ram}.

\ignore{ 
For Ramsey numbers involving cycles, Bondy and Erd\H{o}s \cite{be} showed that $r(C_n, C_n) = 2n - 1$ for odd $n\geq 5$ and $r(C_n, C_{2r-1}) = 2n - 1$ if $n>r(2r-1)$. 
For Ramsey numbers involving paths, Faudree, Lawrence, Parsons and Schelp \cite{flps} showed that $r(P_n, C_m) = 2n - 1$ for $n \ge m \ge 3$ and $r(P_n, C_m) = \max \{ 2n - 1, m + \lfloor n/2 \rfloor - 1 \}$ for $m \ge n \ge 2$ and $m$ is odd.

For Ramsey numbers involving cycles of length at least $m$, Faudree, Schelp and Simonovits \cite{fss} proved that $r(T_n, C_{\ge m}) \le 2m + 2n - 7$ for all $m,n \ge 3$, $r(T_n, C_{\ge m}) \le m + n - 2$ if either $m \ge n$ or $n \ge 432 m^6 - m^2$, and $r(T_n, C_{\ge m}) = n + \lfloor m/2 \rfloor - 1$ if $T_n$ is a tree with maximum degree less than $n - 3m^2$ and $n \ge 432 m^6$. 
Since $P_m$ is a subgraph of $C_{\geq m}$, the upper bounds on $r(T_n, C_{\ge m})$ imply those on $r(T_n, P_m)$ under the same restrictions.
See \cite{rad} for more Ramsey numbers involving paths and cycles.
}


Let $K_k$, $T_k$, $P_k$, $S_k$ and $C_k$ denote the complete graph, tree, path, star and cycle with $k$ vertices, respectively. A classical result is as follows.
\begin{theorem}[Chv\'{a}tal \cite{cha}]\label{chv}
For any integers $k,n\ge1$, $r(T_n,K_k)=(k-1)(n-1)+1$.
\end{theorem}

As pointed out by Burr and Erd\H{o}s \cite{bu-e}: ``Although this theorem is quite easy to prove (in various ways), it seems to 
occupy a central place in Ramsey theory, since it lends itself to many generalizations and analogs." 
For a graph $H$, let  $\chi(H)$ be the chromatic number of $H$, and let $s(H)$ be the minimum size of a color class over all proper vertex-colorings of $H$ by $\chi(H)$ colors. Burr \cite{bur} first observed that if $H_1$ is a graph with chromatic number $\chi(H_1)=k$  and  $H_2$ is a connected graph of order $n\geq s(H_1)$, then
\begin{align}\label{burr}
r(H_1,H_2)\ge(k-1)(n-1)+s(H_1).
\end{align}
Following Burr and Erd\H{o}s \cite{bur,bu-e}, we say a connected graph $H_2$ is {\em $H_1$-good} if the equality of (\ref{burr}) holds. We refer the reader to the survey \cite[Section 2.5]{cfs-15} by Conlon, Fox and Sudakov and \cite{abs,bps,c-l,fl,fox,hhkl,kls,lin,liu,mpy,nr09,ps,ps2} for many recent related references.

Theorem \ref{chv} implies that the tree $T_n$ is $K_k$-good, and it is known \cite{befrs0} that the tree $T_n$  can be replaced by a connected sparse graph $G$ with $n$ vertices.  
For the case $k=3$, Chv\'{a}tal's result implies that $r(T_n,K_3)=2n-1$. Burr et al. \cite[Theorem 14]{befrs1} extended this result by replacing $T_n$ with a connected sparse graph $G$ on $n$ vertices and $K_3$ with an
odd cycle $C_k$ as follows.

\begin{theorem}[Burr, Erd\H{o}s, Faudree, Rousseau and Schelp \cite{befrs1}]\label{Burcycle}
Let $G$ be a connected graph on $n$ vertices, with at most $(1+\frac{1}{42k^5})n$ edges. Then for odd $k\geq 3$ and $n\geq 756k^{10}$, $r(G,C_k)=2n-1$.
\end{theorem}

As a special case, they proposed the following problem.
\begin{problem}[Burr et al. \cite{befrs1}]\label{be-pb}
What is the smallest $n_0$ such that $r(T_n,C_k)=2n-1$  for  odd $k\geq 3$ and $n\ge n_0$?
\end{problem}


Brennan \cite{bre} showed that $n\geq 25k$ will suffice.
\begin{theorem}[Brennan \cite{bre}]
$r(T_n,C_k)=2n-1$  for odd $k\geq 3$ and $n\geq 25k$. 
\end{theorem}

In this paper, we first improve Brennan's bound by reducing the requirement on $n$ from $n \geq 25k$ to $n \geq 4k-8$, with a substantially simpler proof. 
\begin{theorem}\label{main-0}
$r(T_n,C_k)=2n-1$ for odd $k\geq 3$ and $n\geq 4k-8$. 
\end{theorem}

Faudree, Lawrence, Parsons and Schelp \cite{flps} showed that $r(P_n, C_k) = \max \{ 2n - 1, k + \lfloor \frac n2 \rfloor - 1 \}$ for $k \ge n \ge 2$ where $k$ is odd.
In particular, when $n \le 2(k-1)/3$, the second term dominates, giving $r(P_n, C_k) = k + \lfloor n/2 \rfloor - 1 > 2n-1$, which implies $n_0 > 2(k-1)/3$.
Combining this with the upper bound $n_0 \le 4k-8$ given in Theorem \ref{main-0}, we have $$2(k-1)/3 < n_0 \le 4k-8.$$


For the Ramsey goodness problem of sparse graphs versus odd cycles, trees already achieve the optimal linear bound $n = \Omega(k)$ \cite{bre}. In contrast, Theorem \ref{Burcycle} requires $n = \Omega(k^{10})$ for general sparse graphs---a huge polynomial gap that naturally raises the question of what the optimal asymptotic bound should be for sparse graphs in general.

We answer this question completely (up to constant factors), proving that connected sparse graphs are $C_k$-good when $n = \Omega(k)$. This resolves the asymptotic version of the natural extension of Problem \ref{be-pb} to sparse graphs, improving Theorem \ref{Burcycle} in two crucial ways: the vertex requirement drops from $n = \Omega(k^{10})$ to the optimal $n = \Omega(k)$, and simultaneously the edge-density condition relaxes from $1+O(1/k^5)$ to $1+O(1/k^2)$.


 
\begin{theorem}\label{main-2}
Let $G$ be a connected graph on $n$ vertices, with at most  $(1+\frac{1}{36k^2})n$ edges. Then for odd $k\geq 3$ and $n\geq 7k-4$, $r(G,C_k)=2n-1$.
\end{theorem}

As a corollary of Theorem \ref{main-2}, we have the following result. 

\begin{corollary}\label{c-c}
For odd $k\geq 3$ and $n\geq 7k-4$, $r(C_n,C_k)=2n-1$.
\end{corollary}

The classical result, independently due to Bondy and Erd\H{o}s \cite{be}, Faudree and Schelp \cite{f-s}, and Rosta \cite{ros}, states that $r(C_n,C_k) = 2n - 1$ for all $n\ge k\geq 3$ with $k$ odd and $(n,k)\neq(3,3)$ (The condition in \cite{be} is $n> k(k+1)/2$). Consequently, the condition $n\geq 7k-4$ in Corollary \ref{c-c} is tight up to a constant factor.

\medskip
In this paper, we also consider the Ramsey numbers for sparse graphs and paths.
Exact Ramsey numbers are known in several settings involving paths: 
\begin{itemize}
    \item $r(P_n, P_k) = n + \lfloor k/2\rfloor - 1$ for $n \ge k \ge 2$, due to Gerencs\'er and Gyárf\'as~\cite{gg}.
    \item $r(C_n, P_k) = n + \lfloor k/2\rfloor - 1$ for  $k \ge 2$ and $n \ge 3k/2$,  due to Faudree et al.~\cite{flps}.
    \item $r(S_n, P_k) = n + k - 2 - \gamma$ for $k \ge 2$ and $n \ge (k-3)^2+1$, where $\gamma = 0$ if $(k-1) \mid (n-2)$ and $\gamma = 1$ otherwise,  due to Parsons~\cite{par}.
\end{itemize}

The following proposition establishes the tightness of the condition $n \ge (k-3)^2 + 1$ for the equality $r(S_n, P_k) = n + k - 2 - \gamma$. The proof is given in the Appendix.
\begin{proposition} \label{tight-path}
Let $N = n + k - 4$ for $k\ge5$. If the Diophantine equation  
$a(k-1) + b(k-2) = N$ has no non-negative integer solutions $(a,b)$, then $$r(S_n, P_k)\leq N.$$
In particular, when $n = (k-3)^2$  no such solution exists.
\end{proposition}

Burr et al. \cite[Theorem 6]{befrs1}
generalized these results \cite{gg,flps,par} by replacing $P_n,C_n,S_n$ by a sparse graph. 
For any vertex $v$ of a graph $G$, let  $G_v$ be the subgraph of $G$ obtained by deleting $v$ and its neighborhood. Let $\alpha'(G)$ be the minimum independence number among all $G_v$.

\begin{theorem}[Burr, Erd\H{o}s, Faudree, Rousseau and Schelp \cite{befrs1}]\label{Burtree}
Let $G$ be a connected graph on $n$ vertices, with at most $(1+\frac{1}{81k^5})n$ edges. Then for $k\geq 2$ and $n\geq 352k^{12}$, $$r(G,P_k)=\max\left\{n+\lfloor k/2\rfloor-1, n+k-2-\alpha'(G)-\gamma\right\},$$
where $\gamma=0$ if $k-1$ divides $n+k-3-\alpha'(G)$ and $\gamma=1$ otherwise.
\end{theorem}

In this paper, we strengthen Theorem \ref{Burtree} in two ways: we substantially improve the order requirement from $n = \Omega(k^{12})$ to $n = \Omega(k^2)$, while simultaneously relaxing the edge-density condition from $1+O(1/k^5)$ to $1+O(1/k^2)$. In particular, the bound $n = \Omega(k^2)$ is tight up to a constant factor from Proposition \ref{tight-path}.


\begin{theorem}\label{main}
Let $G$ be a connected graph on $n$ vertices, with at most  $(1+\frac{1}{54k^2})n$ edges. Then for $k\geq 2$ and $n\geq 20k^2$, $$r(G,P_k)=\max\left\{n+\lfloor k/2\rfloor-1, n+k-2-\alpha'(G)-\gamma\right\},$$
where $\gamma=0$ if $k-1$ divides $n+k-3-\alpha'(G)$ and $\gamma=1$ otherwise.
\end{theorem}

Since $\alpha'(S_n)=0$, we have the following corollary immediately from the above theorem. 
\begin{corollary}\label{star-path}
For $k \ge 2$ and $n\geq 20k^2$, $r(S_n, P_k) = n + k - 2 - \gamma$, where $\gamma = 0$ if $(k-1) \mid (n-2)$ and $\gamma = 1$ otherwise.
\end{corollary}

If $\alpha'(G)$ is at least $k/2$, then Theorem \ref{main} implies that $r(G,P_k)=n+\lfloor k/2\rfloor-1$.
Let $G$ be a graph satisfying the conditions in Theorem \ref{main}, and $\Delta(G)<n-\frac{\sqrt{n}}{2}$. Then for any vertex $v$ of $G$, $G_v$ has $t> \frac{\sqrt{n}}{2}-1\geq \frac{\sqrt{n}}{4}$ vertices and at most $t+\frac{n}{54k^2}$ edges. Thus $G_v$ has average degree at most $2+\frac{\sqrt{n}}{6k^{2}}.$ Therefore, Tur\'{a}n's Theorem implies that $G_v$ has independence number at least $\frac {t}{1+2+\frac{\sqrt{n}}{6k^{2}}}\geq k/2$, and so $\alpha'(G)\ge k/2$.

Now we obtain a strengthening of Corollary 7 of Burr et al. \cite{befrs1} as follows.

\begin{corollary}\label{genernal}
Let $G$ be a connected graph on $n$ vertices, with at most $(1+\frac{1}{54k^2})n$ edges. If the maximum degree $\Delta(G)< n-\frac{\sqrt{n}}{2}$, $k\geq 2$ and $n\geq 20k^2$, then $$r(G,P_k)=n+\lfloor k/2\rfloor-1.$$
\end{corollary}

If $G$ is a tree on $n$ vertices and $v$ is a vertex of $G$ with degree less than $n-k$, then $G_v$ has $t\geq k$ vertices and at most $t-1$ edges, implying that $G_v$ has independence number at least ${t}/{2}\geq k/2$. As an immediate corollary of Theorem \ref{main}, we obtain a strengthening of Corollary 8 of Burr et al. \cite{befrs1}.

\begin{corollary}\label{tree}
If $T_n$ is a tree with $n$ vertices and $\Delta(T_n)< n-k$, $k\geq 2$ and $n\geq 20k^2$, then $$r(T_n,P_k)=n+\lfloor k/2\rfloor-1.$$
\end{corollary}



In particular, a linear bound on $n$ holds for connected sparse graphs $G$ that satisfy $\delta(G) \ge 2$ and $\alpha'(G) \ge k/2$.
\begin{theorem}\label{main-3}
Let $G$ be a connected graph on $n$ vertices, with minimum degree $\delta(G)\geq2$, $\alpha'(G)\geq \frac{k}{2}$, and at most  $(1+\frac{1}{18k})n$ edges. Then for $k\geq 2$ and $n\geq 3k$, $$r(G,P_k)=n+\lfloor k/2\rfloor-1.$$

\end{theorem}

Since $\alpha'(C_n)\geq \frac{n-3}{2}\geq \frac{k}{2}$ for $n\geq 3k$, we obtain the following result from Theorem \ref{main-3}.
\begin{corollary}\label{p-p-c}
For $k\geq 2$ and $n\geq 3k$, $r(P_n,P_k)=r(C_n,P_k)=n+\left\lfloor\frac k2\right\rfloor-1$.
\end{corollary}

The above corollary implies that $C_n$ ($P_n$) is $P_k$-good for $k\geq 2$ and $n\geq 3k$. Thus the requirement $n\geq 3k$ in Corollary \ref{p-p-c} is tight up to a constant factor.

\medskip
{\em Remark.} We do not attempt to optimize the constants; rather, we emphasize the asymptotic improvement. The constants could likely be improved with more careful analysis.

\medskip

The rest of this paper is organized as follows. In Section \ref{pre}, we present several preliminary results, including a lemma for reconstructing the matching of end-edges (Lemma \ref{findmatch}), which is inspired by the work of Burr et al. \cite{befrs1}, and an enhancement of the dichotomy lemma from the same paper (Lemma \ref{dichotomy}). Building on Lemma \ref{findmatch} and Lemma \ref{dichotomy}, together with other useful tools, we establish in Section \ref{newupper} a new upper bound for $r(G,P_k)$ when $G$ is a connected sparse graph. 
The proofs of Theorem \ref{main-0} (whose proof is much simpler than that of Brennan \cite{bre}) and Theorem \ref{main-2} are provided in Section \ref{pf-0} and Section \ref{pf-1}, respectively, 
and those of Theorem \ref{main} and Theorem \ref{main-3} are given in Section \ref{pf-2}. 
We conclude with remarks and open problems in Section \ref{clu}.

We outline the key ideas leading to the asymptotically optimal bounds in Theorem~\ref{main-2} and Theorem~\ref{main}.

The optimal bound in Theorem~\ref{main-2} is obtained in two stages. First, Lemmas~\ref{sparse-path} and~\ref{sparse-cycle} provide enhanced upper bounds that already apply under the condition $n = \Omega(k^2)$. Second, for the intermediate range $c_1k \le n \le c_2k^2 - 1$ (with $c_1,c_2 > 0$), every connected $n$-vertex graph having at most $(1 + \frac{1}{c_2k^2})n$ edges is either a tree or a tree plus one edge ($T_n + e$). This range is covered by Theorem~\ref{main-0} and Lemma~\ref{special}: Theorem~\ref{main-0} gives a concise alternative proof via an intricate iterative embedding, and Lemma~\ref{special} extends the method to $T_n + e$, showing that $T_n + e$ is $C_k$-good when $n$ is linear in odd $k$. Combining all these elements lowers the vertex requirement from $n = \Omega(k^{10})$ to the optimal $n = \Omega(k)$.

The optimal bound in Theorem~\ref{main} relies on three key innovations: an enhanced upper bound (Lemma~\ref{sparse-path}), a reconstruction technique for end-edge matchings (Lemma~\ref{findmatch}), and a refined dichotomy lemma (Lemma~\ref{dichotomy}). Concretely, whereas Lemma~3 of Burr et al.~\cite{befrs1} demands an end-edge matching of size $\Omega(k^2)$, Lemma~\ref{findmatch} relaxes this to $\Omega(k)$. Simultaneously, Lemma~\ref{dichotomy} strengthens the original dichotomy lemma by raising the lower bound on the number of leaves from $\Omega(n/s)$ to $\Omega(n)$, thereby removing the dependence on $s$ in the leading term. Together, these advances cut the vertex requirement from $\Omega(k^{12})$ to the optimal $\Omega(k^2)$.

\section{Preliminaries}\label{pre}

Let $G=(V,E)$ be a graph. We denote by $|G|$ and $e(G)$ the number of vertices and edges in $G$, respectively. For a subset $A\subseteq V(G)$, let $G[A]$ be the subgraph induced by $A$, and let $G\setminus A$ denote the subgraph induced by $V(G)\setminus A$. If $H$ is a subgraph of $G$, we also write $G \setminus H$ for the subgraph induced by $V(G) \setminus V(H)$.
For a vertex $v\in V(G)$, let $d_G(v)$ be its degree in graph $G$, and let $\delta(G)$ and $\Delta(G)$ denote the minimum and maximum degree of $G$, respectively. For $A \subseteq V(G)$, let $N_G(A)$ denote the set of vertices in $V(G) \setminus A$ that are adjacent to at least one vertex of $A$. For a vertex $v\in V(G)$ and $A \subseteq V(G)$, let $N_A(v)$ denote the set of neighbors of $v$ that lie in $A$.
Let $\alpha(G)$ denote the independence number of $G$. A suspended path is a path whose internal vertices all have degree $2$. An end-edge is an edge incident to a leaf (a vertex of degree $1$). The graph $G + e$ is obtained from $G$ by adding a new edge $e$ between two non-adjacent vertices of $G$. For an edge $e$ of $G$, the graph $G - e$ is obtained from $G$ by deleting $e$.
For $X, Y\subseteq V(G)$, let $E_G(X, Y)$ be the set of edges with one end in $X$ and the other in $Y$. When $X=Y$, we write $E_G(X)$ for $E_G(X, Y)$. A graph $G$ is called pancyclic if it contains cycles of all lengths from 3 to $|V(G)|$.

We will use the following well-known lemma.
\begin{lemma} [Hall \cite{hal}]\label{hall}
A bipartite graph $G:= G(X,Y)$ has a matching which covers every vertex in $X$ if
and only if $$|N_G(S)|\geq|S|\; for\; all \;S\subseteq X.$$
\end{lemma}

Bondy \cite{bon} showed that a graph with sufficiently high minimum degree is pancyclic.
\begin{lemma} [Bondy \cite{bon}]\label{pancyclic}
If an $n$-vertex graph $G$ has $\delta(G)\geq n/2$, then either $G$ is pancyclic or $G$ is the complete bipartite graph $K_{\frac{n}{2}, \frac{n}{2}}$.
\end{lemma}

Observing that the end vertices of a longest path $P$ in a graph $\Gamma$ are completely non-adjacent to the vertices of $V(G\setminus P)$, we have the following result.
\begin{lemma} [Burr et al. \cite{befrs1}]\label{addedge}
For any graph $G$ and for $k\geq 2$, $r(G + e, P_k)\leq r(G, P_k)+k-1$.
\end{lemma}

The following result, which is best possible in some sense, states that if the graph $G$ contains a path whose length cannot be extended by 1, then the complement of $G$ contains a long path. 
\begin{lemma} [Burr et al. \cite{befrs1}]\label{findpath}
Let $u$ and $v$ be distinct vertices of a graph $G$ on $s+t$ vertices where $t\geq s\geq 2$. 
Assume $G$ contains a $t$-vertex path $P$ joining $u$ and $v$ but $G$ contains no $(t+1)$-vertex path joining $u$ and $v$.

\medskip
(i) If $t\geq 2s$, then $\overline{G}$ contains a $P_{2\lceil\frac s2\rceil+1}$.

\medskip
(ii) If $t\geq 5s-1$, then $\overline{G}$ contains a $P_{2s+1}$.

\medskip
(iii) If $t\geq 5s-1$, then $\overline{G}$ contains a $P_{2s-1}$ between any pair of vertices not in $P$ or $\overline{G}$ contains a $K_{2s+1}$.
\end{lemma}

Burr et al. \cite[Lemma 3]{befrs1} proved a lemma for reconstructing the stars of end-edges, thereby requiring the size of the matching of end-edges to be at least $\Omega(k^2)$. We present the following lemma, which reconstructs the end-edge matching directly and relaxes the condition to the optimal $\Omega(k)$.

\begin{lemma}\label{findmatch}
Let $m\geq n\geq 2l\geq 2$ and $k\geq3$ be positive integers. Let $G$ be a graph on $n$ vertices with a matching of $l$ end-edges, and let $H$ be a graph obtained from $G$ by removing $l$ leaves from a matching of $l$ end-edges.
If $r(H, P_k)\leq m-\lfloor\frac{3k}{2}\rfloor$, then $r(G, P_k)\leq m+ \lfloor \frac k2\rfloor-1$. 
\end{lemma}
\noindent{\bf Proof.} 
Let $\Gamma$ be a graph on vertex set $V$ of size $m +\lfloor \frac k2\rfloor-1$. Suppose for contradiction that $\Gamma$ contains no $G$ and $\overline{\Gamma}$ contains no $P_k$. Select a path  $P$ of maximal length in $\overline{\Gamma}$. Delete the vertices of this path and then repeat this procedure. Let $A$ be the set of vertices of the deleted paths. If $|A|\geq k$, we terminate the procedure. Since $|P|\leq k-1$, we have that $k\leq|A|\leq2k-2$ when the procedure finishes.
Since $|\Gamma\setminus A|\geq (m+ \lfloor k/2\rfloor-1)-(2k-2)\geq m-\lfloor\frac{3k}{2}\rfloor\geq r(H, P_k)$, $\Gamma$ contains a copy of $H$ vertex-disjoint from $A$.
By the maximality of the length of $P$, the endvertices of $P$ are adjacent in $\Gamma$ to every vertex not in $A$. 
If $1\leq|P|\leq 2$, then every vertex in $V\setminus A$ is completely adjacent in $\Gamma$ to every vertex on $P$.
If $3\leq|P|\leq k-1$, then by noting for any two consecutive vertices $v_1$, $v_2$ on $P$, every vertex in $V\setminus A$ is adjacent to $v_1$ or $v_2$, we obtain that every vertex in $V\setminus A$ has at least $2+(|P|-3)/2\geq|P|/2$ neighbors on $P$.
Therefore, we conclude that every vertex in $V\setminus A$ has at least $|A|/2$ neighbors in $A$.

Let $L=\{v_1, v_2, \ldots, v_l\}$ be the set of vertices in $H$ that are adjacent to the removed $l$ leaves.
If there are $l$ distinct vertices $\{w_1, w_2, \ldots, w_l\}$ in $V\setminus V(H)$ such that $w_i$ is adjacent in $\Gamma$ to $v_i$ for $1\leq i\leq l$, then $\Gamma$ contains $G$ as a subgraph. If not, by Lemma \ref{hall} we obtain that there exists a set $S =\{v_{i_1}, v_{i_2}, \ldots, v_{i_s}\}\subseteq L$ such that $$|N_\Gamma(S)\cap (V\setminus V(H))|<s.$$
Note that $S\subseteq V(H)\subseteq V\setminus A$ and $A\subseteq V\setminus V(H)$, we obtain that $$s>|N_\Gamma(S)\cap (V\setminus V(H))|\geq |A|/2\geq k/2.$$
Since $|N_\Gamma(S)\cup V(H)|<s+(n-l)\leq n$, we have $|V\setminus (N_\Gamma(S)\cup V(H))|\geq \lfloor k/2\rfloor$. 
Note that each vertex of $S$ is completely non-adjacent to $V\setminus (N_\Gamma(S)\cup V(H))$ in $\Gamma$, we obtain a path $P_k$ in $\overline{\Gamma}$, leading to a contradiction.
\hfill$\Box$

\medskip
Burr et al. \cite[Lemma 2]{befrs1} obtained a dichotomy lemma stating that if a graph $G$ has no long suspended path, then $G$ must have many leaves, i.e., vertices of degree $1$. Combining with the idea of Zhang and Chen \cite[Lemma 4]{zc}, we have a refined dichotomy lemma as follows.
Specifically, in the original proof, Burr et al. replace each suspended path of $G$ with an edge and then later reconstruct each edge back into a suspended path.
We introduce an additional step that deletes all leaves of $G$ at the beginning.
As a result, the bound in the conclusion of the dichotomy lemma is improved from $\Omega(\frac{n}{s})$ to the optimal $\Omega(n)$, removing the dependency on $s$ in the leading term.
\begin{lemma}\label{dichotomy} 
Let $s\ge2$, $n\geq s+1$ and $t\geq -1$ be integers. Let $G$ be a connected graph on $n$ vertices and $n+t$ edges which contains no suspended path with more than $s$ vertices, and let $p$ be the total number of vertices that are adjacent to some leaves in $G$. Then 
$G$ has at least $\min\{n-1, n-(s-1)(2p+3t)-1\}$ leaves. In particular, $G$ has at least $\frac{n-3(s-1)t-1}{2s-1}$ leaves.
\end{lemma}
\noindent{\bf Proof.} 
If $G$ is a star, then $G$ has $n-1$ leaves. 
We may assume that $G$ is not a star in the following.
Let $A$ be the set of vertices that are adjacent to some leaves in $G$. Then $|A|=p$. Let $B\subseteq V(G)\setminus A$ be the set of vertices with degree at least $3$ in $G$, and $C\subseteq V(G)\setminus A$ the set of all vertices of degree $2$  in $G$.
Let $G_1$ be the graph obtained from $G$ by removing all leaves. 
Let $G_2$ be the graph (possibly a multigraph) obtained from $G_1$ by applying the following two operations:
\medskip

(i) if $G_1$ contains a suspended path with all internal vertices from $C$ and endvertices $u, v$ from $A\cup B$, then we remove these internal vertices and add a new edge between $u$ and $v$;
\medskip

(ii) if $G_1$ contains a cycle with exactly one vertex $v$ from $A\cup B$ and all other vertices from $C$, then we remove all vertices of the cycle except $v$ and add a loop to $v$.
\medskip

If $G$ contains a cycle in which all vertices belong to $C$, then $G$ itself is a cycle since $G$ is connected. This contradicts the assumption that $G$ has no suspended path with more than $s$ vertices. Thus any cycle in $G$ has at least one vertex that does not belong to $C$.

From the above two operations, we remove all vertices in $C$ from $G_1$.
Thus 
$$e(G_1)=|G_1|+t=|A|+|B|+|C|+t, \;\;\text{and}\;\; e(G_2)=|G_2|+t=|A|+|B|+t.$$
Since $G$ is connected and not a star, every vertex of $A$ has degree at least $1$ in $G_1$. Also, every vertex of $B$ has degree at least $3$ in $G_1$, and every vertex of $C$ has degree $2$ in $G_1$. Thus $|A|+3|B|+2|C|\leq 2(|G_1|+t)=2(|A|+|B|+|C|+t)$, and so we obtain $|B|\leq |A|+2t$.

Note that $G_1$ can be reconstructed from $G_2$ by replacing each non-loop edge with a suspended path of length at most $s-1$ (keeping the original endpoints) and each loop with an induced cycle of length at most $s$.
During the reconstruction, for each non-loop edge in $G_2$, we add at most $s-2$ vertices; for each loop in $G_2$, we add at most $s-1$ vertices. Since $G_2$ is connected, $G_2$ has at least $|G_2|-1$ non-loop edges and at most $e(G_2)-|G_2|+1$ loops. It follows that
$$|G_1|\le |G_2|+(s-2)(|G_2|-1)+(s-1)(e(G_2)-|G_2|+1)= (s-1)e(G_2)+1.$$
Thus the total number of leaves in $G$ is 
\begin{align*}
n-|G_1|&\geq n-(s-1)e(G_2)-1\\
&=n-(s-1)(|A|+|B|+t)-1\\
&\geq n-(s-1)(2|A|+3t)-1\\
&= n-(s-1)(2p+3t)-1,
\end{align*}
where the last inequality follows by noting $|B|\leq |A|+2t$.
Moreover, by the definition of $p$, $G$ has at least $p$ leaves. Thus when $n-(s-1)(2p+3t)-1=p$, i.e., when $p=\frac{n-3(s-1)t-1}{2s-1}$, the number of leaves in $G$ attains the minimum value $\frac{n-3(s-1)t-1}{2s-1}$. 

Combining the two cases, we obtain that $G$ has at least $\min\{n-1, n-(s-1)(2p+3t)-1\}$ leaves. In particular, $G$ has at least $\frac{n-3(s-1)t-1}{2s-1}$ leaves.
\hfill$\Box$

\section{A new upper bound on $r(G,P_k)$}\label{newupper}
 

We employ an upper bound on $r(T_n,P_k)$ established by Häggkvist \cite{hag}, which will be used in the base case of the induction for the upper bound of $r(G, P_k)$ in Lemma \ref{sparse-path}.

\begin{lemma}[Häggkvist \cite{hag}]\label{tree-path}
$r(K_{n,m}, P_k)\leq n+m+k-2$.
\end{lemma}

Since every tree is a subgraph of a complete bipartite graph, the following corollary holds.
\begin{corollary}\label{base}
For integers $k\geq 2$ and $n\geq2$, $r(T_n, P_k)\leq  n+k-2$.
\end{corollary}

The above upper bound is also valid for $T_n + e$.
\begin{lemma}\label{unicyclic-path}
For integers $k\geq2$ and $n\geq 3k-2$, $r(T_n+e, P_k)\leq n+k-2$.
\end{lemma}
\noindent{\bf Proof.} 
Let $T_n$ be a tree on $n$ vertices with bipartition $(A,B)$, and let $G = T_n + e$. Then $|A|+|B|=n$. Consider a graph $\Gamma$ on $n+k-2$ vertices, and suppose for contradiction that $\Gamma$ does not contain $G$ and $\overline{\Gamma}$ does not contain $P_k$.
By Lemma \ref{tree-path}, $\Gamma$ contains $K_{|A|,|B|}$ as a subgraph; denote its partite sets by $(U,W)$, where $|U| = |A|$ and $|W| = |B|$.

We now examine the position of $e$ in $G$.  
If $e \in E_G(A,B)$, then $\Gamma$ already contains a copy of $G$, a contradiction.  
If $e \in E_G(A)$ and $|U| \ge k$, then either an edge in $\Gamma[U]$ allows us to embed $G$ into $\Gamma$, or $\overline{\Gamma}[U]$ is a clique of order at least $k$ and hence contains a $P_k$; in either case we obtain a contradiction.  
The case $e \in E_G(B)$ with $|W| \ge k$ is symmetric, yielding a contradiction as well.
Thus we may assume that $G$ contains an odd cycle $C$.  

Let $x$ be an endvertex of $e$ in $G$, and let $T(x)$ be the component containing $x$ of the subgraph induced by $V(G)\setminus(V(C)\setminus\{x\})$.
By symmetry, assume $e \in E_G(A)$. Thus $|A|=|U|\leq k-1$ from the above discussion.

If $|W| = |B| \ge |V(G \setminus T(x)) \cap B| \ge k$, we swap the vertices in $A\cap V(T(x))$ with $B\cap V(T(x))$.
Let $e'$ be the other edge of $C$ incident to $x$. Let $T_n' = G - e'$ with two parts $A', B'$. Then the vertices of $V(G \setminus T(x)) \cap B$ still lie in $B'$, i.e., $B'\supseteq V(G \setminus T(x)) \cap B$.
By Lemma \ref{tree-path}, $\Gamma$ contains $K_{|A'|,|B'|}$ as a subgraph; denote its partite sets by $(U',W')$, where $|U'| = |A'|$ and $|W'| = |B'|$.    
After this reassignment, we have $e' \in E_G(B')$ and $|W'| = |B'| \ge |V(G \setminus T(x)) \cap B| \ge k$, which yields a contradiction by a similar argument as above.  
Therefore
\[
|V(G \setminus T(x)) \cap B| \le k-1.
\]

If $|W| = |B| \ge |V(T(x)) \cap B| \ge k$, we swap $A$ with $B$ for all vertices in $G\setminus T(x)\cup\{x\}$ and keep all other vertices in $T(x)\setminus \{x\}$ fixed.
Similarly, we obtain that $$|V(T(x)) \cap B| \le k-1.$$ 
Therefore, we have
$
|V(G)| = |A|+|V(G \setminus T(x)) \cap B|+|V(T(x)) \cap B|\le 3k-3.
$
This contradicts the assumption that $|V(G)| = n \ge 3k-2$.
\hfill$\Box$

\medskip
We are now in a position to prove Lemma \ref{sparse-path}, which improves the upper bound of Burr et al. \cite[Proposition 5]{befrs1}.
\begin{lemma}\label{sparse-path}
For integers $k\geq 3$ and $n\geq 3k$, let $G$ be a connected graph on $n$ vertices, with at most $(1+\frac{1}{27k^2})n$ edges. Then $r(G, P_k)\leq n+k-2$.
\end{lemma}
\noindent{\bf Proof.} 
The proof will be by induction on $n$. 
If $3k\leq n\leq 27k^2-1$, then $$e(G)\leq \left(1+\frac{1}{27k^2}\right)n\leq n+\frac{27k^2-1}{27k^2}<n+1,$$
implying that $G$ has at most $n$ edges. 
Since $G$ is connected, $G$ must be either a tree $T_n$ or $T_n+e$. 
By Corollary \ref{base} and Lemma \ref{unicyclic-path}, the lemma follows.
Thus we may assume $n\geq27k^2$ and the result is true for all $m$-vertex appropriate graphs with $3k\leq m\leq n-1$. Let $\Gamma$ be a graph on $n+k-2$ vertices. Suppose for contradiction that $\Gamma\nsupseteq G$ and $\overline{\Gamma}\nsupseteq P_k$. 

\medskip\noindent
{\bf Case 1} \; $G$ contains a suspended path with $3k$ vertices.

\medskip
Let $H$ be the graph on $n-k$ vertices obtained from $G$ by shortening the suspended path by $k$ vertices. Thus the number of edges of $H$ is
at most $$\left(1+\frac{1}{27k^2}\right)n-k\leq \left(1+\frac{1}{27k^2}\right)(n-k)+1.$$
Since $r(H-e, P_k)\leq n-k+k-2= n-2$ by induction, applying Lemma \ref{addedge} gives $r(H, P_k)\leq n+k-2$. Thus $\Gamma$ contains $H$ as a subgraph.
Let $H^+$ be a subgraph of $\Gamma$ in which this suspended path has been lengthened as much as possible (may be up to $k$ vertices). If $H^+=G$, then we are done, otherwise select a vertex set $S$ of size $k-1$ of $\Gamma$ not in $H^+$. Let $\Gamma'$ be the subgraph of $\Gamma$ induced by the vertices of such lengthened suspended path of $H^+$ and $S$. Applying Lemma \ref{findpath}(i) with $t\geq 3k-k=2k\geq2|S|$, we obtain that $\overline{\Gamma'}$, and hence $\overline{\Gamma}$, contains a $P_{2\lceil\frac {k-1}2\rceil+1}\supseteq P_k$, which leads to a contradiction.

\medskip\noindent
{\bf Case 2} \; Each suspended path of $G$ has at most $3k-1$ vertices.

\medskip
Let $p$ be the total number of vertices that are adjacent to some leaves in $G$.
Suppose first that $p\geq 3k$, then $G$ contains a matching of $3k$ end-edges. 
Let $H$ be the graph on $n-3k$ vertices obtained from $G$ by removing $3k$ leaves from the matching of $3k$ end-edges. Thus the number of edges of $H$ is at most $(1+\frac{1}{27k^2})n-3k\leq (1+\frac{1}{27k^2})(n-3k)+1$. 
By the induction assumption, we have $r(H-e, P_k)\leq n-3k+k-2=n-2k-2$. From Lemma \ref{addedge}, we obtain 
$r(H, P_k)\leq n-k-2$.
Applying Lemma \ref{findmatch} with $m=n+\lfloor\frac{k}2\rfloor-2$, we obtain that $r(G, P_k)\leq n+k-2$ as desired. 
So we may assume that $p< 3k$ in the following.

Applying Lemma \ref{dichotomy} with $s=3k-1$ and $t\leq\frac{n}{27k^2}$, we obtain that $G$ has at least $\min\{n-1, n-(3k-2)(2p+3t)-1\}$ leaves.
Thus $G$ contains a star with at least 
\begin{align*}
\frac{\min\{n-1, n-(3k-2)(2p+3t)-1\}}p&\geq\frac{n-3k(2\cdot3k+3\cdot\frac{n}{27k^2})}{3k}\\
&\geq\frac{n}{3k}\left(1-\frac{1}{3k}\right)-6k\\
&\geq2k
\end{align*} 
end-edges provided $n\geq 27k^2$ and $k\geq 3$.
Let $S$ be a set of $2k$ leaves in $G$ which are adjacent to a vertex $v$ in $G$. Let $u$ be the endvertex of a path $P$ of maximal length in $\overline{\Gamma}$. Thus the path $P$ has at most $k-1$ vertices. Let $H=G\setminus S$. Then $H$ has $n-2k$ vertices and at most $(1+\frac{1}{27k^2})n-2k\leq (1+\frac{1}{27k^2})(n-2k)+1$ edges.
Since $r(H-e, P_k)\leq n-2k+k-2=n-k-2$ by induction, applying Lemma \ref{addedge} gives $r(H, P_k)\leq n-2$. Thus $\Gamma$ contains $H$ as a subgraph such that $H$ is vertex-disjoint from $P$. Since $u$ is adjacent in $\Gamma$ to each vertex of $\Gamma$ not on $P$ and $$d_{\Gamma\setminus(H\setminus\{v\}\cup P)}(u)\geq n+k-2-(n-2k-1)-(k-1)= 2k,$$ replacing $v$ in $H$ by $u$ gives a copy of $G$ in $\Gamma$, which leads to a contradiction.
\hfill$\Box$

\section{Proof of Theorem \ref{main-0}}\label{pf-0}

From (\ref{burr}), we have $r(T_n, C_k)\geq 2n-1$  for each odd $k\geq 3$. In the following, we will show that if $k\geq 3$ is odd and $n\geq 4k-8$, then $r(T_n, C_k)\leq 2n-1$.
Let $\Gamma$ be a graph on a vertex set $V$ of size $2n-1$. Suppose for contradiction that $\Gamma$ contains no $T_n$ and $\overline{\Gamma}$ contains no $C_k$.

If $\delta(\Gamma) \ge n - 1$, then $T_n$ can be embedded into $\Gamma$. Thus $\delta(\Gamma) \le n - 2$, which implies that $\Delta(\overline{\Gamma}) = 2n - 2 - \delta(\Gamma) \ge n$. Let $v$ be a vertex of $\Gamma$ with $d_{\overline{\Gamma}}(v) = \Delta(\overline{\Gamma}) \ge n$, and 
let $S\subseteq N_{\overline{\Gamma}}(v)$ with $|S|=n$.
Since $\overline{\Gamma}$ is $C_k$-free, $\overline{\Gamma}[S]$ is $P_{k-1}$-free.
If $k=3$, then $\Gamma[S]$ forms a complete graph on $n$ vertices, which is impossible.
So we may assume $k\geq5$. 
Let $H_1$ be a connected subtree obtained from $T_n$ by repeatedly deleting leaves until exactly $n-k+3$ vertices remain.
By Corollary \ref{base}, $r(H_1, P_{k-1})\leq (n-k+3)+(k-1)-2= n.$
Thus $\Gamma[S]$ contains a copy of $H_1$, which we will denote by $F_1$. 

We define a greedy procedure to extend this embedding of $H_1$ to an embedding of $T_n$ in $\Gamma$. During this procedure, we maintain a subgraph $K_1\subseteq \Gamma$ induced by the set of vertices that have already been embedded. Initially, we have $K_1 = F_1$.
The embedding is extended iteratively according to the following rule: 
Whenever we have an edge $x_1x_2 \in E(T_n)$ and a vertex $v \in V$ such that $x_1$ is already mapped to $v \in V(K_1)$ but $x_2$ is not yet embedded, we map $x_2$ to some vertex in $N_{V(\Gamma \setminus K_1)}(v)$, provided this neighborhood is nonempty.
Note that each $x \in V(T_n)$ that is not yet embedded has at most one neighbor among the already embedded vertices of $T_n$. 
Since $\Gamma$ does not contain $T_n$ as a subgraph, this embedding procedure must fail. 
Consequently, $|V(K_1)| \le n-1$ and there exists a vertex $v_1\in V(K_1)$ with $N_{V(\Gamma \setminus K_1)}(v_1) = \emptyset$.
Let $T=V\setminus S$ and $U_1=V(\Gamma\setminus K_1)$. Then $|U_1|\geq n$. Moreover, $|V(K_1)\cap T|\leq k-4$ since otherwise we can embed $T_n$ into $\Gamma$. Thus
$$|U_1\cap T|\geq (n-1)-(k-4)=n-k+3.$$

We will select inductively $\frac{k-1}{2}$ vertices $v_1, v_2, \ldots, v_{\frac{k-1}{2}}$ from $\Gamma$ such that the corresponding set $U_i$ $(1\leq i\leq\frac{k-1}{2})$ consisting of vertices non-adjacent to $v_i\in \Gamma \setminus \{v_1,\ldots, v_{i-1}\}$ satisfies that $|U_i|\geq n-(i-1)$ and $|U_i\cap T|\geq n-k+5-2i$. The base case holds from the above.
Suppose that $\{v_1,\ldots, v_{i-1}\}$ have been selected for $2\leq i\leq \frac{k-1}{2}$. We now select $v_i$ from $\Gamma\setminus \{v_1,\ldots, v_{i-1}\}$.
For $2\leq i\leq \frac{k-1}{2}$, let $H_i$ be a connected subtree of $T_n$ on $n-k+4-i$ vertices, and let $S_i=S\setminus \{v_1,\ldots, v_{i-1}\}$.
By Corollary \ref{base}, $r(H_i, P_{k-1})\leq n-k+4-i+k-3\leq n-(i-1).$
Thus $\Gamma[S_i]$ contains a copy of $H_i$, which we will denote by $F_i$. 
We use the above greedy procedure to extend this embedding of $H_i$ to an embedding of $T_n$ in $\Gamma\setminus \{v_1,\ldots, v_{i-1}\}$. 
During this procedure, we maintain a subgraph $K_i\subseteq \Gamma$ induced by the set of vertices that have already been embedded. Initially, we have $K_i = F_i$. Similarly, we conclude that $|V(K_i)| \le n-1$ and there exists a vertex $v_i\in V(K_i)$ with $N_{U_i}(v_i) = \emptyset$ where $$U_i=V \setminus (\{v_1,\ldots, v_{i-1}\}\cup V(K_i)).$$ 
Then $|U_i|\geq n-(i-1)$.
Note that $$|(\{v_1,\ldots, v_{i-1}\}\cup V(K_i))\cap T|\le(n-1)-(n-k+4-i)+(i-1)=k+2i-6.$$ Hence $|U_i\cap T|\geq n-1-(k+2i-6)=n-k+5-2i$.

Let $W=\{v_1, v_{\frac{k-1}{2}}, v_3, v_{\frac{k-5}{2}}, \ldots, v_{\frac{k-7}{2}},v_4, v_{\frac{k-3}{2}},v_2\}$, 
and let $p_i$ be the subscript of the $i$-th vertex in the sequence $W$ for $1\leq i\leq \frac{k-1}{2}$.
Then $p_i+p_{i+1}\leq\frac{k-1}{2}+3$, and
\begin{align*}
|U_{p_i}\cap U_{p_{i+1}}\cap T|&\geq n-1-(k-5+p_i)-(k-5+p_{i+1})-\max\{p_i-1,p_{i+1}-1\}\\
&\geq n-2k+9-\left(\frac{k-1}{2}+3\right)-\frac{k-3}{2}\\
&\geq k
\end{align*}
provided $n\geq 4k-8$.
For any vertices $u_1\in U_1$ and $u_2\in U_2$, if $u_1u_2\in E(\overline{\Gamma})$, then $\overline{\Gamma}$ contains a 
$$C_k=u_1v_1w_1v_{\frac{k-1}{2}}w_2v_3w_3v_{\frac{k-5}{2}}\cdots v_{\frac{k-7}{2}}w_{\frac{k-7}{2}}v_4w_{\frac{k-5}{2}}v_{\frac{k-3}{2}}w_{\frac{k-3}{2}}v_2u_2u_1$$ 
where each $w_i$  for $1\leq i\leq \frac{k-3}{2}$  is chosen from $U_{p_i} \cap U_{p_{i+1}} \cap T$ and the order in which all vertices $v_i$ of $C_k$ appear is given by the sequence $W$.
Thus we can assume that $u_1u_2\in E(\Gamma)$. It follows that $U_1\cap U_2$ induces a complete graph in $\Gamma$ and each vertex of $U_1$ is adjacent in $\Gamma$ to each vertex in $U_1\cap U_2$. 
In particular, $|U_1\cap U_2|\geq n-2k+5$ and $|U_1|\geq n$.
Consider a bipartition $V(T_n)=A\cup B$ where $|A|\ge |B|$. It follows that $|A|\ge n/2$.
Now consider the embedding of $T_n$ into $\Gamma[U_1]$. Note that $|U_1\setminus U_2| \le 2k-5 \le n/2 \le |A|$ since $n \ge 4k-10$.  We can injectively embed the $|U_1\setminus U_2|$ vertices in $A$ into $U_1\setminus U_2$ and embed the remaining $|U_1\cap U_2|$ vertices of $T_n$ into $U_1\cap U_2$. 
This contradicts the fact that $\Gamma$ does not contain $T_n$ as a subgraph. The proof of Theorem \ref{main-0} is complete.
\hfill$\Box$

\section{Proof of Theorem \ref{main-2}}\label{pf-1}

We first prove the following result, which is of independent interest and plays a crucial role in the proof of Theorem \ref{main-2}.
\begin{lemma}\label{special}
$r(T_n+e,C_k)=2n-1$ for odd $k\geq 3$ and $n\geq 7k-4$. 
\end{lemma}
\noindent{\bf Proof.} 
It suffices to prove the upper bound, because the lower bound is already given by (\ref{burr}).
Let $G = T_n + e$ with $e = xy$, and let $(A, B)$ be a bipartition of $T_n$.
Note that $G$ contains exactly one cycle, say $C$, and let $m = |C|$ be its length, where $3 \le m \le n$.
Let $\Gamma$ be a graph on a vertex set $V$ of size $2n-1$. 
Suppose for contradiction that $\Gamma$ does not contain $G$ as a subgraph and $\overline{\Gamma}$ does not contain $C_k$ as a subgraph.

We claim that $\Delta(\overline{\Gamma})\geq n-1$. Otherwise, we have $\delta(\Gamma)\geq |V|-(n-1)=n>|V|/2$.
From Lemma \ref{pancyclic}, $\Gamma$ contains $C$ as a subgraph.
We apply the greedy procedure from the proof of Theorem \ref{main-0} to extend this embedding of $C$ to an embedding of $G$ in $\Gamma$.
During this procedure, we maintain a subgraph $K\subseteq \Gamma$ induced by the set of vertices that have already been embedded.
Initially, we have $K = C$.
Similarly, we conclude that $|V(K)| \le n-1$ and there exists a vertex $v_0\in V(K)$ with $N_{V(\Gamma\setminus K)}(v_0) = \emptyset$. 
This implies that $d_{\overline{\Gamma}}(v_0)\geq n$, a contradiction.

Let $u_0\in V$ be a vertex  of maximum degree in $\overline{\Gamma}$.
Let $S$ be a subset of $N_{\overline{\Gamma}}(u_0)$ with $n-1$ vertices, and let $T=V\setminus S$.
Since $\overline{\Gamma}$ does not contain $C_k$, $\overline{\Gamma}[S]$ contains no $P_{k-1}$.
If $k=3$, then $\Gamma[S]$ would be a complete graph on $n$ vertices, a contradiction.
So we may assume $k\geq5$. 

Suppose that $m\leq n-\frac{3k-7}{2}$. Let $U_1,U_2$ be defined as in Theorem \ref{main-0}.
Then by an argument analogous to that in Theorem \ref{main-0}, we can apply Lemma \ref{unicyclic-path} to embed $G$ into $\Gamma[U_1]$, provided $n\geq \frac{9(k-1)}{2}$. This contradicts the fact that $\Gamma$ contains no $G$.
The main differences are the following. First, the graph $H_i$ we choose is a connected subgraph of $G$ on $n-k+3-i$ vertices containing $C$ for $1\leq i\leq \frac{k-1}{2}$. Second, when embedding $G$ into $\Gamma[U_1]$, we injectively embed the $|U_1\setminus U_2|$ vertices of $A\setminus \{x,y\}$ into $U_1\setminus U_2$.

Next we assume that $m\geq n-\frac{3(k-3)}{2}$.
Note that $|A\cap C|, |B\cap C|\geq \frac{m-1}{2}$. Since each vertex in $V(G)\setminus V(C)$ has at most one neighbor in $C$,
without loss of generality, we may assume that at most half vertices of $V(G)\setminus V(C)$ are adjacent to some vertices in $B\cap C$.
Let $I$ be a subset of $B\cap C\setminus(\{x,y\}\cup N_C(\{x,y\}))$ such that
each vertex of $I$ has exactly degree two in $G$ and each pair of vertices in $I$ has no neighbors in common.
Note that there are at least $|B\cap C|-4-\frac{n-m}{2}$ vertices in $B\cap C\setminus\{x,y\}$ which have degree two in $G$ from the above assumption.
Thus $$|I|\geq \left\lfloor\frac{1}{2}\left(|B\cap C|-4-\frac{n-m}{2}\right)\right\rfloor\geq \frac{n-3(k-3)-9}{4}-1\geq k-2$$ provided $n\geq 7k-4$.
Let $I'\subseteq I$ with $|I'|=k-2$. 

We now embed $G\setminus I'$ into $\Gamma[S]$ at first.
Since $\overline{\Gamma}[S]$ contains no $P_{k-1}$, from Lemma \ref{tree-path}, $\Gamma[S]$ contains $K_{|A|,|B|-|I'|}$ as a subgraph; denote its partite sets by $(U,W)$, where $|U| = |A|$ and $|W| = |B|-|I'|$.
Thus we can embed $G\setminus I'$ into $\Gamma[U\cup W]$, except that the edge $e$ remains to be mapped.
We now examine the position of $e$ in $G$.  
If $e \in E_G(A,B\setminus I')$, then $\Gamma[S]$ already contains a copy of $G\setminus I'$.  
Provided $m\geq n-\frac{3(k-3)}{2}$ and $n\geq\frac{11k-19}{2}$, we obtain $|W| \geq \frac{m-1}{2}-(k-2)\geq k-1$ and $|U| \geq \frac{m-1}{2}\geq k-1$.
Suppose $e \in E_G(B \setminus I')$. 
If there is an edge in $\Gamma[W]$, we can then embed $G \setminus I'$ into $\Gamma[S]$. Otherwise, $\overline{\Gamma}[W]$ is a clique of order at least $k-1$, then $\overline{\Gamma}[W]$ contains a $P_{k-1}$, leading to a contradiction.
Suppose  $e \in E_G(A)$. Then we can  embed $G \setminus I'$ into $\Gamma[S]$ similarly.
In all cases, $\Gamma[S]$ contains a copy $H$ of $G\setminus I'$, which we will denote by $F$. 

We define a greedy procedure to extend this embedding of $G\setminus I'$ to an embedding of $G$ in $\Gamma$. 
During this procedure, we maintain a subgraph $K\subseteq \Gamma$ induced by the set of vertices that have already been embedded. Initially, we have $K = F$.

The embedding is extended iteratively using the following rule:
Whenever we have two edges $x_1x_3, x_2x_3\in E(G)$ and two vertices $v_1, v_2\in V$ such that $x_1, x_2$ are already mapped to $v_1, v_2\in V(K)$ respectively but $x_3$ is not yet embedded, we map $x_3$ to some vertex in $N_{V(\Gamma \setminus K)}(v_1)\cap N_{V(\Gamma \setminus K)}(v_2)$, provided this neighborhood is nonempty. 
Note that each $x \in I'$ has exactly two neighbors among the already embedded vertices of $G$. 
Let $u_1, u_2, \ldots, u_t$, with $t\leq k-3$, be the images in $\Gamma$ of the vertices of $I'$ that have already been embedded in this procedure,
and let $u, v$ be the images in $\Gamma$ of $x,y\in V(G)$ under the same embedding map.
Define $U'=U\setminus (\cup_{i=1}^t N_U(u_i)\cup \{u,v\})$.
Then $|U'|\geq \frac{m-1}{2}-2(k-3)-2$.
Since $\Gamma$ does not contain $G$ as a subgraph, this embedding procedure must fail. 
Consequently, we have that $|V(K)| \le n-1$ and any vertices in $V\setminus V(K)$ have at most one neighbor in $U'$ by noting that $\Gamma(U,W)$ is a complete bipartite graph.
If $V\setminus V(K)$ forms a complete graph of $|V\setminus V(K)|\ge n$ vertices, then there is a copy of $G$, leading to a contradiction.
Thus there exists a non-edge $v_1 v_2$ in $\Gamma\setminus K$.
For any two vertices $v'$, $v''$ in $V\setminus V(K)$, 
$$|N_{\overline{\Gamma}[U']}(v')\cap N_{\overline{\Gamma}[U']}(v'')|\geq |U'|-2\geq\frac{m-1}{2}-2(k-3)-2-2\geq \frac{k-1}{2}$$ 
provided $m\geq n-\frac{3(k-3)}{2}$ and $n\geq \frac{13k-17}{2}$.
Therefore, we can greedily find a $C_k$ using $\{v_1, v_2\}$, $\frac{k-1}{2}$ vertices from $U'$ and $\frac{k-3}{2}$ vertices from $V\setminus (V(K)\cup\{v_1, v_2\})$,  a contradiction.
\hfill$\Box$

\medskip
We use the following result to bound Ramsey numbers of general graphs versus cycles.
\begin{lemma} [Burr, Erd\H{o}s, Faudree, Rousseau and Schelp \cite{befrs1}]\label{upbound}
If $G$ is a graph on $n$ vertices and $l$ edges, then $r(G, C_k)\leq n+2lk-\frac{2l}{n}$.
\end{lemma}

Before proving Theorem \ref{main-2}, we first derive a new upper bound on $r(G,C_k)$.

\begin{lemma}\label{sparse-cycle}
For any integers $k\geq 3$ and $n\geq 3k+1$, if $G$ is a connected graph with $n$ vertices and at most $(1+\frac{1}{27k^2})n$ edges, then $r(G, C_k)\leq 2n+k-2$.
\end{lemma}
\noindent{\bf Proof.} 
Assume the assertion is not true. Then there exists a graph $\Gamma$ on a vertex set $V$ of size $2n + k-2$ such that $\Gamma \nsupseteq G$ and $\overline{\Gamma} \nsupseteq C_k$. Let $G'$ be the graph obtained from $G$ by successively deleting
leaves and shortening suspended paths with at least $3k+1$ vertices by a vertex. 
 \begin{claim} 
$\Gamma$ does not contain $G'$ as a subgraph.
\end{claim}

\noindent{\bf Proof.}  We prove the assertion by repeatedly applying the following two cases.

\medskip\noindent
{\bf Case 1} \; $G$ has a leaf.

\medskip
Let $H$ be the graph obtained from $G$ by deleting a leaf. Let $v$ be the vertex of $G$ adjacent to this leaf. If $\Gamma$ contains $H$ as a subgraph, then the vertex $v$ is adjacent in $\overline{\Gamma}$ to at least $n+k-1$ vertices. By Lemma \ref{sparse-path}, we obtain that $r(G, P_k)\leq n+k-2$. Thus the neighborhood of $v$ in $\overline{\Gamma}$ contains  $P_k$, which together with $v$ yields a $C_k$ in $\overline{\Gamma}$.
So $H$ is not a subgraph of $\Gamma$.

\medskip\noindent
{\bf Case 2} \; $G$ contains a suspended path $P$ with at least $3k+1$ vertices.

\medskip
Let $H$ be the graph obtained from $G$ by shortening the suspended path $P$ by a vertex. Let $s=\lfloor\frac{k+1}{2}\rfloor$. Since $2s+1\ge k$ and $\overline{\Gamma}$ contains no $C_k$, we have that $\overline{\Gamma}$ contains no $K_{2s+1}$. Suppose that $\Gamma$ contains $H$.
Thus $\Gamma$ contains a path of length at least $3k$ whose length cannot be extended by 1. 
We apply Lemma \ref{findpath}(iii) to $\Gamma$ with $t\geq 3k\geq 5s-1$ to obtain that for any pair of vertices in $V\setminus V(H)$, 
$\overline{\Gamma}$ contains either a path with $k$ vertices between them if $k$ is odd or $k-1$ vertices between them otherwise.
By Lemma \ref{sparse-path}, we obtain that $r(G, P_k)\leq n+k-2\leq |V\setminus V(H)|$ where $k\geq 3$. 
Thus $\overline{\Gamma}$ contains a $P_3=(v_1, v_2, v_3)$ vertex-disjoint from $H$.
Let $k'=\lfloor\frac{k+1}{2}\rfloor-2$ and let $w_1, w_2, \ldots,w_{k'}$ be any $k'$ vertices in $V\setminus (V(H\cup P_3))$.

If $k$ is odd, let $\Gamma'$ be the subgraph of $\Gamma$ spanned by the vertices of the suspended path of $H$ and $\{v_1, v_2, w_1, w_2, \ldots,w_{k'}\}$. From Lemma \ref{findpath}(iii), $\overline{\Gamma}$ contains a path with $k$ vertices between $v_1$ and $v_2$. Thus $\overline{\Gamma}$ contains a $C_k$, leading to a contradiction.

If $k$ is even, let $\Gamma'$ be the subgraph of $\Gamma$ spanned by the vertices of the suspended path of $H$ and $\{v_1, v_3, w_1, w_2, \ldots,w_{k'}\}$. From Lemma \ref{findpath}(iii), $\overline{\Gamma}$ contains a path with $k-1$ vertices between $v_1$ and $v_3$, which together with $P_3$ forms a $C_k$ in $\overline{\Gamma}$, again a contradiction.
We can thus conclude that $H$ is not a subgraph of $\Gamma$. \hfill$\Box$

\medskip
Suppose that the graph $G'$ has $l$ vertices. Then $G'$ contains at most $l+\frac{n}{27k^2}$ edges, and it has no suspended path with more than $3k$ vertices, and no leaves.

Applying Lemma \ref{dichotomy} with $s=3k$ and $t\leq\frac{n}{27k^2}$, and using the fact that $n\geq 3k+1\geq s+1$,  we obtain that $G'$ has at least 
$\frac{l-3(s-1)t-1}{2s-1}$ leaves. Since $G'$ has no leaves, it follows by noting 
$k\geq3$ and  $\frac{n}{3k}\geq1$ that
$l\leq3(s-1)t+1\leq 3\cdot3k\cdot \frac{n}{27k^2}+\frac{n}{3k}=\frac{2n}{3k}.$
Thus, by Lemma \ref{upbound}, 
$r(G',C_k)\leq \frac{2n}{3k}+2k\left(\frac{2n}{3k}+\frac{n}{27k^2}\right)\leq 2n.$ Therefore, $\Gamma$ contains $G'$ as a subgraph, contradicting the above claim.
\hfill$\Box$

\medskip
Now we are ready to prove Theorem \ref{main-2}.

\medskip
\noindent
{\bf Proof of Theorem \ref{main-2}.}\; 
From (\ref{burr}), we obtain $r(G, C_k)\geq 2n-1$. In the following, we will show that for any connected graph $G$ with $n$ vertices and at most $(1+\frac{1}{36k^2})n$ edges, provided $k\geq 3$ is odd  and $n\geq 7k-4$, then $r(G, C_k)\leq 2n-1$.
Let $\Gamma$ be a graph on vertex set $V$ of size $2n-1$. Suppose for contradiction that $\Gamma\nsupseteq G$ and $\overline{\Gamma}\nsupseteq C_k$.
If $7k-4\leq n\leq 18k^2-1$, then
$$e(G)\leq \left(1+\frac{1}{36k^2}\right)n\leq n+\frac{18k^2-1}{36k^2}<n+1,$$
implying that $G$ has at most $n$ edges. 
Since $G$ is connected, $G$ must be either a tree $T_n$ or $T_n+e$. 
By Theorem \ref{main-0} and Lemma \ref{special}, the result follows.
Thus we may assume $n\geq18k^2$. 

\medskip\noindent
{\bf Case 1} \; $G$ contains a suspended path with at least $3k+1$ vertices.

\medskip
Let $H$ be the graph on $n-\frac{k-1}{2}$ vertices obtained from $G$ by shortening the suspended path by $\frac{k-1}{2}$ vertices. Thus $H$ has at most $(1+\frac{1}{36k^2})n-\frac{k-1}{2}\leq (1+\frac{1}{27k^2})(n-\frac{k-1}{2})$ 
edges provided $n\geq 2k-2$.
It follows by Lemma \ref{sparse-cycle} that $r(H, C_k)\leq 2(n-\frac{k-1}{2})+k-2= 2n-1$, which implies that $\Gamma$ contains $H$ as a subgraph. Let $H^+$ be a subgraph of $\Gamma$ in which this suspended path has been lengthened as much as possible (up to $\frac{k-1}{2}$ vertices). If $H^+=G$, then we are done. If not, applying Lemma \ref{findpath}(iii) with $t\geq \frac{5(k+1)}{2}-1:=5s-1$, we obtain that for any pair of vertices in $V\setminus V(H)$, $\overline{\Gamma}$ contains a path $P_k$ between them. Since $\overline{\Gamma}$ does not contain a $C_k$, any pair of vertices in $V\setminus V(H)$ must be adjacent in $\Gamma$. This yields a complete graph on at least $n$ vertices in $\Gamma$ and hence $G$ is a subgraph of $\Gamma$, leading to a contradiction.

\medskip\noindent
{\bf Case 2} \; Each suspended path of $G$ has at most $3k$ vertices.

\medskip
Applying Lemma \ref{dichotomy} with $s=3k$ and $t\leq\frac{n}{36k^2}$, we obtain that $G$ has at least $$\frac{n-3(s-1)t-1}{2s-1}\geq \frac{n-3\cdot3k\cdot\frac{n}{36k^2}-1}{6k}\geq\frac{n}{6k}-\frac{n}{24k^2}-1\geq3k-3$$ leaves provided $n\geq 18k^2$ and $k\geq 3$. 

Let $H$ be the graph obtained from $G$ by deleting $\frac{3(k-1)}{2}$ leaves. 
Thus the number of edges of $H$ is at most $(1+\frac{1}{36k^2})n-\frac{3(k-1)}{2}\leq (1+\frac{1}{27k^2})(n-\frac{3(k-1)}{2})$ provided $n\geq 6k-6$.
It follows by Lemma \ref{sparse-cycle} that $r(H, C_k)\leq2(n-\frac{3(k-1)}{2})+k-2\leq 2n-1.$ Thus $\Gamma$ contains $H$ as a subgraph. Let $H^+$ be a maximal connected subgraph of $G$ containing $H$ which is a subgraph of $\Gamma$. If $H^+ = G$, we are done. Otherwise, there is a vertex $v$ of $H^+$ which is adjacent in $\overline{{\Gamma}}$ to each vertex in $V\setminus V(H^+)$. Let $S$ be a subset of $V\setminus V(H^+)$ with $n$ vertices and let $$T = V\setminus S.$$

Since $\overline{\Gamma}$ does not contain $C_k$, $\overline{\Gamma}[S]$ contains no $P_{k-1}$.
If $k=3$, then $\Gamma[S]$ forms a complete graph on $n$ vertices. Thus $\Gamma$ contains $G$ as a subgraph, a contradiction.
So we may assume that $k\geq5$. Since $e(H)\leq (1+\frac{1}{27k^2})(n-\frac{3(k-1)}{2})$, by Lemma \ref{sparse-path}, we obtain that 
\begin{align}\label{H-Pk-1}
r(H, P_{k-1})\leq n-\frac{3(k-1)}{2}+k-2\leq n-\frac{k-1}{2}.
\end{align} 
Thus $\Gamma[S]$ contains a copy of $H$, which we will denote by $F_1$. 
By adding leaves, we enlarge $F_1$ as much as possible to obtain a subgraph $F_1^+$ of $\Gamma$. 

If $F_1^+= G$,
we are done. Otherwise, there is a vertex $v_1$ of $F_1$ which is adjacent in $\overline{\Gamma}$ to each vertex in $V\setminus V(F_1^+)$. 
Let $U_1=V\setminus V(F_1^+)$. Hence $|U_1|\geq n$ and $$|U_1\cap T|\geq (n-1)-\left(\frac{3(k-1)}{2}-1\right)=n-\frac{3(k-1)}{2}.$$ 

From (\ref{H-Pk-1}), we can select inductively $\frac{k-1}{2}$ vertices $v_1, v_2, \ldots, v_{\frac{k-1}{2}}$ from $S$ such that the corresponding set $U_i$ $(1\leq i\leq\frac{k-1}{2})$ consisting of vertices adjacent to $v_i$ in $\overline{\Gamma}$ satisfies that $|U_i|\geq n$ and $|U_i\cap T|\geq n-\frac{3(k-1)}{2}$.

Let $U=\cap_{i=1}^{(k-1)/2} (U_i\cap T)$. Then $U$ has at least $n-\frac{k-1}{2}\cdot\frac{3(k-1)}{2}\ge k$ vertices. Therefore there exists a complete bipartite graph in $\overline{\Gamma}$ with two parts $\{v_1, v_2, \ldots, v_{\frac{k-1}{2}}\}$ and  $U$. For any vertices $u_1\in U_1$ and $u_2\in U_2$, if $u_1u_2\in E(\overline{\Gamma})$, then $\overline{\Gamma}$ contains a $C_k$ using $\{u_1, u_2, v_1, v_2, \ldots, v_{\frac{k-1}{2}}\}$ and $\{w_1, w_2, \ldots, w_{\frac{k-3}{2}}\}$ from $U$. 
Thus we can assume that $u_1u_2\in E(\Gamma)$. It follows that $U_1\cap U_2$ induces a complete graph  in $\Gamma$ and each vertex of $U_1$ is adjacent in $\Gamma$ to each vertex in $U_1\cap U_2$. Since $|U_1\cap U_2|\geq n-3(k-1)$ and
$|U_1|\geq n$, we can embed $G$ into $\Gamma$ by noting $G$ has at least $3k-3$ leaves. This completes the proof of Theorem \ref{main-2}.
\hfill$\Box$

\medskip
{\em Remark.} In the context of Ramsey numbers for cycles versus sparse graphs of minimum degree at least two, a slightly better constant factor can be obtained.

\begin{proposition}\label{main-4}
Let $G$ be a connected graph on $n$ vertices, with minimum degree $\delta(G)\geq2$ and at most  $(1+\frac{1}{36k^2})n$ edges. Then for odd $k\geq 3$ and $n\geq 3k+1$, $r(G,C_k)=2n-1$.
\end{proposition}
\noindent{\bf Proof.} 
Similar to the proof of Lemma \ref{sparse-path}, we obtain that for integers $k\geq 3$ and $n\geq 3k$, if $G$ is a connected graph on $n$ vertices with minimum degree $\delta(G)\geq2$ and at most $(1+\frac{1}{9k})n$ edges, then $r(G, P_k)\leq n+k-2$. The main difference lies in the argument from Case 2 in the proof of Lemma \ref{sparse-path}. Applying Lemma \ref{dichotomy} with $p=0$, $s=3k-1$ and $t\leq\frac{n}{9k}$, we obtain that $G$ has at least $\min\{n-1, n-3(s-1)t-1\}\geq n-3(3k-2)\cdot\frac{n}{9k}-\frac{n}{3k}=\frac{n}{3k}>0$ leaves, which contradicts the assumption that $\delta(G)\geq2$.
Then similar to the proof of Lemma \ref{sparse-cycle}, we obtain that for integers $k\geq 3$ and $n\geq 3k+1$, if $G$ is a connected graph on $n$ vertices with minimum degree $\delta(G)\geq2$ and at most $(1+\frac{1}{18k^2})n$ edges, then $r(G, C_k)\leq 2n+k-2$. 

Now we are ready to prove Proposition \ref{main-4}.
For the lower bound, we have $r(G,C_k)\ge 2n-1$ from (\ref{burr}).
For the upper bound, by an argument similar to that in Case 1 of the proof of Theorem \ref{main-2}, we conclude that each suspended path of $G$ has at most $3k$ vertices.
Applying Lemma \ref{dichotomy} with $p=0$, $s=3k$ and $t\leq\frac{n}{36k^2}$, we obtain that $G$ has at least $\min\{n-1, n-3(s-1)t-1\}\geq n-3\cdot3k\cdot\frac{n}{36k^2}-\frac{n}{3k}\geq(1-\frac{7}{12k})n>0$ leaves, which contradicts the assumption that $\delta(G)\geq2$.
\hfill$\Box$

\medskip
Proposition \ref{main-4} yields a slight strengthening of Corollary \ref{c-c}.
\begin{corollary}
For odd $k\geq 3$ and $n\geq 3k+1$, $r(C_n,C_k)=2n-1$.
\end{corollary}

\section{Proofs of Theorem \ref{main} and Theorem \ref{main-3}}\label{pf-2}

In this section, we shall prove Theorem \ref{main} and Theorem \ref{main-3}.

\begin{proposition}\label{struc-f} 
If $\Gamma$ contains no $P_k$, then either $\Gamma$ is a vertex-disjoint union of connected graphs, each with at most $k-1$ vertices, or $\Gamma$ contains a vertex of degree at most $k/2-1$.
\end{proposition}
{\bf Proof.} Consider a maximal length path $P =(x_1, x_2,\ldots, x_t)$ in $\Gamma$. Then $t\leq k-1$ since there is no $P_k$ in $\Gamma$.
If there exists an $i$ $(1 < i\leq t)$ such that $x_1x_i$ and $x_{i-1}x_t$ are edges in $\Gamma$, then the vertices of $P$ form a cycle $(x_1, x_i, x_{i+1},\ldots, x_t, x_{i-1}, x_{i-2}, \ldots, x_2, x_1)$.
The maximality of the length of $P$ implies that the vertices of $P$ form a connected component of $\Gamma$. If no such $i$ exists, then
for each $2\leq i\leq t-1$, either $x_1x_i\in {\overline{\Gamma}}$ or $x_{i-1}x_t\in {\overline{\Gamma}}$. Since $x_1x_2$ and $x_{t-1}x_t$ are edges in $\Gamma$, we obtain $d_{\Gamma}(x_1)+d_{\Gamma}(x_t)\leq 2+(t-3)= t-1\leq k-2$. Therefore, either $x_1$ or $x_t$ has degree at most $k/2-1$ in $\Gamma$. The proposition follows. 
\hfill$\Box$
\medskip

\noindent{\bf Proof of Theorem \ref{main}.} 
The proof of the lower bound is similar as in \cite{befrs1,par}. We include the proof for completeness.
From (\ref{burr}),  $r(G, P_k)\geq n+\left\lfloor\frac k2\right\rfloor-1$. 
Next, if $\alpha'(G)<k/2$, we show that $r(G, P_k)\geq n+k-2-\alpha'(G)-\gamma$, where $\alpha'(G)$ is the minimum independence number among all subgraphs $G_v$ for $v\in V(G)$, and $\gamma=0$ if $k-1$ divides $n+k-3-\alpha'(G)$ and $\gamma=1$ otherwise.
From \cite[Lemma 5]{par}, we obtain that if $N\geq (k-2)^2-(k-2)$, then there exist non-negative integers $p$, $q$ such that $N=p(k-2)+ q(k-1)$.
Let $N=n+k-3-\alpha'(G)-\gamma$. Then $N\geq (k-2)^2-(k-2)$ provided $n\geq20k^2\geq (k-3)^2+k/2+1$. 
Let $\Gamma$ be a complete $(p+q)$-partite graph where the first $p$ parts have size $k-2$ and the remaining $q$ parts have size $k-1$.
Thus $\overline{\Gamma}=pK_{k-2}\cup qK_{k-1}$. In particular, $\overline{\Gamma}=qK_{k-1}$ if $k-1$ divides $n+k-3-\alpha'(G)$.
Since every connected component of $\overline{\Gamma}$ has fewer than $k$ vertices, $\overline{\Gamma}$ contains no $P_k$. We will show that $\Gamma$ does not contain $G$ as a subgraph. Suppose for contradiction that $G$ is a subgraph of $\Gamma$. 
Let $v$ be a vertex in $G$ such that $\alpha'(G)=\alpha(G_v)$. The vertex $v$ lies in one of the disjoint independent sets of $\Gamma$, say $X$. Then $X$ has at least $k-1-\gamma$ vertices by noting $\gamma=0$ if $k-1$ divides $n+k-3-\alpha'(G)$ by the definition of $\gamma$. 
Moreover, we have $|X\cap V(G)|\le \alpha(G_v)+1=\alpha'(G)+1$.
However, then $|\Gamma|=|\Gamma\backslash X|+|X|\geq|(\Gamma\backslash X)\cap V(G)|+|X|\geq n-(\alpha'(G)+1)+(k-1-\gamma)=n+k-2-\alpha'(G)-\gamma>|\Gamma|,$
 leading to a contradiction. 
The proof of the lower bound for $r(G, P_k)$ is complete.

\medskip
It remains to show the upper bound of $r(G, P_k)$.
To this end, let $\Gamma$ be a graph on $N=\max\{n+\left\lfloor\frac k2\right\rfloor-1, n+k-2-\alpha'(G)-\gamma\}$ vertices.
Suppose for contradiction that $\Gamma$ contains no $G$ and $\overline{\Gamma}$ contains no $P_k$.
We assume $k\geq3$ since the assertion is clear for $k=2$. 

\medskip\noindent
{\bf Case 1} \; $G$ contains a suspended path with $3k$ vertices.

\medskip
Let $H$ be the graph on $n-\lfloor\frac {k}2\rfloor$ vertices obtained from $G$ by shortening the suspended path by $\lfloor\frac {k}2\rfloor$ vertices. 
Thus $H$ has at most $(1+\frac{1}{54k^2})n-\lfloor\frac {k}2\rfloor\leq (1+\frac{1}{27k^2})(n-\lfloor\frac {k}2\rfloor)$ edges provided $n\geq k$.
By Lemma \ref{sparse-path}, we obtain that $r(H, P_k)\leq n+\left\lfloor\frac k2\right\rfloor-1\leq N$ which implies that $\Gamma$ contains $H$ as a subgraph.
Let $H^+$ be a subgraph of $\Gamma$ in which the suspended path has been lengthened as much as possible (up to $\lfloor\frac {k}2\rfloor$ vertices). If $H^+=G$, then we are done. If not, applying Lemma \ref{findpath}(ii) with $t\geq\frac {5k}2-1$, we obtain that $\overline{\Gamma}$ contains a $P_{2\lfloor\frac {k}2\rfloor+1}\supseteq P_k$, a contradiction.

\medskip\noindent
{\bf Case 2} \; Each suspended path of $G$ has at most $3k-1$ vertices. 

\medskip
Let $p$ be the total number of vertices that are adjacent to some leaves in $G$.
If $p\geq \frac{5k}{2}$, then $G$ contains a matching of $b$ end-edges where $b=\lfloor\frac{5k}{2}\rfloor$. 
Let $H$ be the graph on $n-b$ vertices obtained from $G$ by removing $b$ leaves from such a matching of $b$ end-edges. Thus $H$ has at most $(1+\frac{1}{54k^2})n-b\leq (1+\frac{1}{27k^2})(n-b)$ edges. 
By Lemma \ref{sparse-path}, $r(H, P_k)\leq n-b+k-2\leq n-\lfloor\frac{3k}{2}\rfloor$.
Applying Lemma \ref{findmatch} with $m=n$, $r(G, P_k)\leq n+\lfloor\frac k2\rfloor-1\leq N$. This completes the proof in this case.
So we may assume that $p< \frac{5k}{2}$ in the following. 

Applying Lemma \ref{dichotomy} with $s=3k-1$ and $t\leq\frac{n}{54k^2}$, 
we obtain that $G$ has at least $\min\{n-1, n-(3k-2)(2p+3t)-1\}$ leaves.
Thus $G$ contains a star with at least 
\begin{align*}
\frac{\min\{n-1, n-(3k-2)(2p+3t)-1\}}{p}\geq\frac{n-3k(2\cdot\frac{5k}{2}+3\cdot\frac{n}{54k^2})}{\frac{5k}{2}}
\geq\frac{2n}{5k}-\frac{n}{15k^2}-6k
\geq\frac{3k}{2}
\end{align*}
end-edges provided $n\geq 20k^2$ and $k\geq3$. 

Let $S$ be a set of $\lfloor\frac{3k}{2}\rfloor$ leaves which are adjacent in $G$ to a vertex $v$ and let $H=G\setminus S$. Thus $H$ has $n-2k$ vertices and at most $(1+\frac{1}{54k^2})n-\lfloor\frac{3k}{2}\rfloor\leq (1+\frac{1}{27k^2})(n-\lfloor\frac{3k}{2}\rfloor)$ edges. Therefore, by Lemma \ref{sparse-path}, $r(H, P_k)\leq n-\lfloor\frac{3k}{2}\rfloor+k-2= n-\lfloor\frac{k}{2}\rfloor-2$.

\medskip
Let $u$ be a vertex of $\overline{\Gamma}$ of minimum degree. If $d_{\overline{\Gamma}}(u)\leq k/2-1$, then $$d_{\Gamma}(u)\geq \left(n+\left\lfloor k/2\right\rfloor-1\right)-\left(\left\lfloor k/2\right\rfloor-1\right)-1=n-1\geq r(H, P_k).$$ Thus there is a copy of $H$ in $N_\Gamma(u)$. Replacing $v$ in $H$ with $u$ gives another copy of $H$ in $\Gamma$. Since $d_{\Gamma}(u)\geq n-1$ and $H=G\setminus S$, we have that $G$ is a subgraph of $\Gamma$, a contradiction.

Thus we assume $d_{\overline{\Gamma}}(u)>k/2-1$. 
From Proposition \ref{struc-f}, $\overline{\Gamma}$ is a vertex-disjoint union of connected graphs, each with at most $k-1$ vertices. 
Select a connected component $W$ of $\overline{\Gamma}$ of minimum size. We claim that if $\Gamma$ has exactly $n+k-2-\alpha'(G)-\gamma$ vertices, then $|W| \leq k-1-\gamma$. This follows because, by assumption, $\gamma=1$ precisely when $k-1$ does not divide $n+k-3-\alpha'(G)$.

Since $r(H, P_k)\leq n-\lfloor\frac{k}{2}\rfloor-2$, 
$\Gamma$ contains a copy $F$ of $H$ which is vertex-disjoint from $W$. 
The edges between $W$ and $F$ are all in $\Gamma$. Let $w$ be any vertex in $W$. Replacing $v$ in $F$ with $w$ gives a copy of $H$ in $\Gamma$, which we will denote by $F_1$. 
Select an independent set $I$ of $\min\{\alpha(F_v), |W| -1\}$ vertices in $F_v$.
Replacing $I$ in $F_1$ with $I_W\subseteq W\setminus\{w\}$ such that $|I_W|=|I|$ yields another copy of $H$, which we will denote by $F_2$. 
Note that each vertex of $W$ is adjacent in $\Gamma$ to each vertex of $F_1$ except $w$. 
The graphs $F$, $F_1$ and $F_2$ can be seen in Fig. 1.

\begin{figure}[t]
\begin{center}

\tikzset{every picture/.style={line width=0.75pt}} 

\begin{tikzpicture}[x=0.75pt,y=0.75pt,yscale=-1,xscale=1]

\draw   (185,145.5) .. controls (185,113.47) and (205.82,87.5) .. (231.5,87.5) .. controls (257.18,87.5) and (278,113.47) .. (278,145.5) .. controls (278,177.53) and (257.18,203.5) .. (231.5,203.5) .. controls (205.82,203.5) and (185,177.53) .. (185,145.5) -- cycle ;
\draw  [fill={rgb, 255:red, 0; green, 0; blue, 0 }  ,fill opacity=1 ] (232,121) .. controls (232,119.62) and (233.12,118.5) .. (234.5,118.5) .. controls (235.88,118.5) and (237,119.62) .. (237,121) .. controls (237,122.38) and (235.88,123.5) .. (234.5,123.5) .. controls (233.12,123.5) and (232,122.38) .. (232,121) -- cycle ;
\draw   (351,176.5) .. controls (351,162.97) and (362.42,152) .. (376.5,152) .. controls (390.58,152) and (402,162.97) .. (402,176.5) .. controls (402,190.03) and (390.58,201) .. (376.5,201) .. controls (362.42,201) and (351,190.03) .. (351,176.5) -- cycle ;
\draw    (263.5,134.5) -- (302.53,145.65) -- (319.5,150.5) ;
\draw   (308,135.5) .. controls (308,95.18) and (335.76,62.5) .. (370,62.5) .. controls (404.24,62.5) and (432,95.18) .. (432,135.5) .. controls (432,175.82) and (404.24,208.5) .. (370,208.5) .. controls (335.76,208.5) and (308,175.82) .. (308,135.5) -- cycle ;
\draw   (336,167.25) .. controls (336,145.85) and (353.24,128.5) .. (374.5,128.5) .. controls (395.76,128.5) and (413,145.85) .. (413,167.25) .. controls (413,188.65) and (395.76,206) .. (374.5,206) .. controls (353.24,206) and (336,188.65) .. (336,167.25) -- cycle ;
\draw  [fill={rgb, 255:red, 0; green, 0; blue, 0 }  ,fill opacity=1 ] (340,113.25) .. controls (340,111.87) and (341.12,110.75) .. (342.5,110.75) .. controls (343.88,110.75) and (345,111.87) .. (345,113.25) .. controls (345,114.63) and (343.88,115.75) .. (342.5,115.75) .. controls (341.12,115.75) and (340,114.63) .. (340,113.25) -- cycle ;
\draw    (345,113.25) -- (371.5,110.5) ;
\draw    (342.5,110.75) -- (362.5,97.75) ;
\draw   (357.5,95.75) .. controls (357.5,84.29) and (369.25,75) .. (383.75,75) .. controls (398.25,75) and (410,84.29) .. (410,95.75) .. controls (410,107.21) and (398.25,116.5) .. (383.75,116.5) .. controls (369.25,116.5) and (357.5,107.21) .. (357.5,95.75) -- cycle ;
\draw    (264,148.5) -- (319,136.5) ;
\draw   (210,167.25) .. controls (210,153.3) and (221.42,142) .. (235.5,142) .. controls (249.58,142) and (261,153.3) .. (261,167.25) .. controls (261,181.2) and (249.58,192.5) .. (235.5,192.5) .. controls (221.42,192.5) and (210,181.2) .. (210,167.25) -- cycle ;
\draw    (378.5,138.5) -- (381.5,110.5) ;
\draw  [dash pattern={on 0.84pt off 2.51pt}]  (392.5,140.5) -- (390.5,112) ;
\draw    (244.07,112.68) .. controls (269.38,84.01) and (299.74,85.42) .. (332.02,107.48) ;
\draw [shift={(333.5,108.5)}, rotate = 215.05] [color={rgb, 255:red, 0; green, 0; blue, 0 }  ][line width=0.75]    (10.93,-3.29) .. controls (6.95,-1.4) and (3.31,-0.3) .. (0,0) .. controls (3.31,0.3) and (6.95,1.4) .. (10.93,3.29)   ;
\draw [shift={(242.5,114.5)}, rotate = 309.99] [color={rgb, 255:red, 0; green, 0; blue, 0 }  ][line width=0.75]    (10.93,-3.29) .. controls (6.95,-1.4) and (3.31,-0.3) .. (0,0) .. controls (3.31,0.3) and (6.95,1.4) .. (10.93,3.29)   ;
\draw    (264.75,181.54) -- (273.3,181.7) -- (341.75,182.96) ;
\draw [shift={(343.75,183)}, rotate = 181.06] [color={rgb, 255:red, 0; green, 0; blue, 0 }  ][line width=0.75]    (10.93,-3.29) .. controls (6.95,-1.4) and (3.31,-0.3) .. (0,0) .. controls (3.31,0.3) and (6.95,1.4) .. (10.93,3.29)   ;
\draw [shift={(262.75,181.5)}, rotate = 1.06] [color={rgb, 255:red, 0; green, 0; blue, 0 }  ][line width=0.75]    (10.93,-3.29) .. controls (6.95,-1.4) and (3.31,-0.3) .. (0,0) .. controls (3.31,0.3) and (6.95,1.4) .. (10.93,3.29)   ;

\draw (188.5,75.9) node [anchor=north west][inner sep=0.75pt]    {$W$};
\draw (220.5,101.9) node [anchor=north west][inner sep=0.75pt]    {$w$};
\draw (371.5,168.9) node [anchor=north west][inner sep=0.75pt]    {$I$};
\draw (227,158.4) node [anchor=north west][inner sep=0.75pt]    {$I_{W}$};
\draw (420.5,62.4) node [anchor=north west][inner sep=0.75pt]    {$F$};
\draw (358,133.9) node [anchor=north west][inner sep=0.75pt]    {$F_{v}$};
\draw (335,91.9) node [anchor=north west][inner sep=0.75pt]    {$v$};
\draw (362,86.9) node [anchor=north west][inner sep=0.75pt]    {$N_{F}( v)$};
\draw (276.5,68.4) node [anchor=north west][inner sep=0.75pt]    {$F_1$};
\draw (286.5,184.4) node [anchor=north west][inner sep=0.75pt]    {$F_2$};

\end{tikzpicture}

\end{center}
\begin{center}
Fig. 1. The graphs $F$, $F_1$ and $F_2$ where $V(F_1) =V(F)\setminus \{v\}\cup \{w\}, V(F_2) =V(F_1)\setminus I\cup I_{W}$.
\end{center}
\end{figure}

If $\alpha(F_v)\geq |W| -1$, then $|I_W|=|I|=|W|-1$ and so $w$ is adjacent  in $\Gamma$ to all vertices in $\Gamma\setminus F_2$. 
Since $d_{\Gamma\setminus F_2}(w)\geq \left(n+\left\lfloor  k/2\right\rfloor-1\right)-(n-\lfloor\frac{3k}{2}\rfloor)\geq\lfloor\frac{3k}{2}\rfloor,$ we obtain that $\Gamma$ contains $G$ as a subgraph, leading to a contradiction. 

Thus we may assume that $\alpha(F_v)< |W| -1$.
Since $\alpha(F_v)=\alpha(H_v)=\alpha(G_v)\geq\alpha'(G)$, where the equality holds as $H$ is obtained from $G$ by deleting some neighbors of $v$, we have that $w$ is adjacent in $\Gamma$ to all except at most $|W|-1-\alpha'(G)$ vertices in $\Gamma\setminus F_2$. 

Note that $N\geq n+k-2-\alpha'(G)-\gamma$. If $N=n+k-2-\alpha'(G)-\gamma$, then
\begin{align*}
d_{\Gamma\setminus F_2}(w)\geq N-\left(n-\left\lfloor\frac{3k}{2}\right\rfloor\right)-(|W|-1-\alpha'(G))\geq \left\lfloor\frac{3k}{2}\right\rfloor
\end{align*}
by noting $|W|\le k-1-\gamma$.
If $N\geq n+k-1-\alpha'(G)-\gamma$, then
\begin{align*}
d_{\Gamma\setminus F_2}(w)\geq N-\left(n-\left\lfloor\frac{3k}{2}\right\rfloor\right)-(|W|-1-\alpha'(G))\geq \left\lfloor\frac{3k}{2}\right\rfloor
\end{align*}
by noting $|W|\le k-1$.
Thus $\Gamma$ contains $G$ as a subgraph, leading to a contradiction. 
This final contradiction completes the proof of Theorem \ref{main}.
\hfill$\Box$

\medskip
\noindent{\bf Proof of Theorem \ref{main-3}.}  For the lower bound,  we have $r(G,P_k)\ge n+\lfloor\frac k2\rfloor-1$ from (\ref{burr}) directly; hence we do not need the requirement $n\ge (k-3)^2+\frac k2+1$ in the proof of the lower bound of Theorem \ref{main}.

For the upper bound, recall that for integers $k\geq 3$ and $n\geq 3k$, if $G$ is a connected graph on $n$ vertices with minimum degree $\delta(G)\geq2$ and at most $(1+\frac{1}{9k})n$ edges, then $r(G, P_k)\leq n+k-2$. 
Then similar to the proof of Theorem \ref{main}, we can finish the proof of Theorem \ref{main-3}.
The main difference lies in the argument from Case 2 in the proof of Theorem \ref{main}. Applying Lemma \ref{dichotomy} with $p=0$, $s=3k-1$ and $t\leq\frac{n}{18k}$, we obtain that $G$ has at least $$\min\{n-1, n-3(s-1)t-1\}\geq n-3\cdot3k\cdot\frac{n}{18k}-1=\frac{n}{2}-1>0$$ leaves, which contradicts the assumption that $\delta(G)\geq2$.
\hfill$\Box$

\section{Concluding remarks}\label{clu}

Our work provides a comprehensive picture of Ramsey goodness for sparse graphs against odd cycles and paths, establishing asymptotically optimal vertex requirements. For odd cycles, Theorem \ref{main-2} determines the exact Ramsey number under near-linear edge density when $n = \Omega(k)$. For paths, Theorem \ref{main} achieves the same under the condition $n = \Omega(k^2)$. Theorem \ref{main-3} shows that a linear bound $n = \Omega(k)$ suffices for graphs with minimum degree at least two and $\alpha'(G) \ge k/2$---a condition that is tight up to a constant factor. The constants appearing in our theorems are not optimized and could likely be improved with more careful analysis, though this would not affect the asymptotic statements.

\medskip
Several natural questions remain:

\medskip
{\bf Problem 1. The exact threshold.} Determine the minimal $n_0(k)$ for which the assertions of Theorems \ref{main-2}, \ref{main} and \ref{main-3} hold. In particular, close the constant-factor gap between the current lower and upper bounds, or prove that $\lim_{k \to \infty} n_0(k)/k$ exists and compute its value.

{\bf Problem 2. Relaxing the edge density.} Can the edge-density assumption be relaxed from $1 + O(1/k^2)$ to $1 + O(1/k)$? Our methods do not achieve this, raising the question of whether the quadratic dependence is intrinsic when graphs $G$ are allowed to have leaves.

{\bf Problem 3. Even cycles--a striking contrast.} The behavior for even cycles is fundamentally different. Even the simplest case $r(S_n, C_4)$ remains an open challenge. Burr, Erd\H{o}s, Faudree, Rousseau, and Schelp \cite{befrs} conjectured (with a \$100 prize offered by Erd\H{o}s) that for any constant $c > 0$, $r(S_{n+1}, C_4) < n + \sqrt{n} - c$ holds for infinitely many $n$. More generally, determining Ramsey numbers for sparse graphs versus even cycles remains wide open, as the corresponding extremal graphs exhibit a delicate and complex structure that resists precise characterization.


\section*{Appendix}

\medskip\noindent
{\bf Proof of Proposition \ref{tight-path}.}
Let $\Gamma$ be a graph on $N = n + k - 4$ vertices.
Suppose for contradiction that $\Gamma$ contains no $S_n$ and $\overline{\Gamma}$ contains no $P_k$.
Let $u$ be a vertex of minimum degree in $\overline{\Gamma}$.  
If $d_{\overline{\Gamma}}(u) \le k/2 - 1$, then we have $d_\Gamma(u) \ge (n + k - 4) - (\lfloor \frac k2\rfloor - 1) - 1 \ge n - 1,$  so $\Gamma$ contains an $S_n$, a contradiction.
Thus we may assume $d_{\overline{\Gamma}}(u) > k/2 - 1$. Then by Proposition \ref{struc-f}, $\overline{\Gamma}$ must be a vertex-disjoint union of connected graphs, each of order at most $k-1$.
Choose a connected component $W$ of $\overline{\Gamma}$ of minimum size.  
Since the Diophantine equation  
$a(k-1) + b(k-2) = N$ has no non-negative integer solutions $(a,b)$, we have $|W| \le k-3$.
Thus any vertex $v$ in $W$ has at most $|W| - 1 \le k-4$ neighbors in $\overline{\Gamma}$, which implies  
$d_\Gamma(v) \ge (n + k - 4) - (k - 3) = n - 1,$ a contradiction.

In particular, for $n=(k-3)^2$, we have $N = n + k - 4 = k^2 - 5k + 5$.
To establish $r(S_n, P_k) \le N$ in this case, it suffices to verify that the Diophantine equation $a(k-1) + b(k-2) = N$ has no non-negative integer solution $(a,b)$. Suppose for contradiction that such a pair $(a,b)$ exists. 
Then $$a+b>a + b\cdot\frac{k-2}{k-1} =\frac{N}{k-1} \;\;\text{and}\;\; a+b<a\cdot\frac{k-1}{k-2} + b =\frac{N}{k-2}.$$
Thus $\frac{N}{k-1}<a+b<\frac{N}{k-2}$.
Note that $N =(k-1)(k-4)+1=(k-2)(k-3)-1,$ 
it follows that $k-3\leq a+b\leq k-4$, leading to a contradiction.
\hfill$\Box$

\end{spacing}

\end{document}